\documentclass{elsarticle}
\usepackage{amsfonts}
\usepackage{amsmath}
\usepackage{mathabx}
\usepackage[all]{xy}
\usepackage{tikz}
\newtheorem{thm}{Theorem}
\newtheorem{lem}[thm]{Lemma}

\newproof{pf}{Proof}
\newproof{pot1}{Proof of Theorem \ref{mainthm}}
\newproof{pot2}{Proof of Theorem \ref{mainthm0}}

\DeclareMathOperator{\interior}{int}

\begin{document}
\begin{frontmatter}

\title{Global convergence in systems of differential equations
arising from chemical reaction networks}

\author[ref1]{Murad Banaji\corref{cor1}} \author[ref2]{Janusz
Mierczy\'nski} \address[ref1]{Department of Mathematics, University
of Portsmouth, Lion Gate Building, Lion Terrace, Portsmouth,
Hampshire PO1 3HF, UK} \address[ref2]{Institute of Mathematics and
Computer Science, Wroc{\l}aw University of Technology, Wybrze\.ze
Wyspia\'nskiego 27, PL-50-370 Wroc{\l}aw, Poland} \cortext[cor1]{{\tt
murad.banaji@port.ac.uk}}

\begin{abstract}
It is shown that certain classes of differential equations arising
from the modelling of chemical reaction networks have the following
property: the state space is foliated by invariant subspaces each of
which contains a unique equilibrium which, in turn, attracts all
initial conditions on the associated subspace.
\end{abstract}

\begin{keyword}
Global attractivity; monotone dynamical systems; chemical reactions;
network structure; DSR graph
\end{keyword}

%\begin{AMS}
%80A30; 34C99; 05C38; 05C50
%\end{AMS}

\end{frontmatter}

\section{Introduction} There are difficult and interesting open
questions about allowed asymptotic behaviour in systems of
differential equations arising in the modelling of chemical reaction
networks (CRNs for short). The main goal in this area is to make
claims about the behaviour of these systems which are as far as
possible independent of the particular choices of functions or
parameters which describe the rates of reaction or ``kinetics''.
Classical results in this direction \cite{hornjackson,feinberg} rely
strongly on the choice of ``mass action kinetics'' leading to
particular polynomial differential equations. Mathematically, such
results involve proving that solutions of certain parameterised
families of polynomial differential equations have certain
asymptotic behaviours regardless of the values of the parameters, but
provided these have fixed sign. However when the kinetics is
constrained only by loose qualitative laws, the family of possible
differential equations describing a reaction network becomes much
larger, and results become fewer. Here, we provide a general result
based on the theory of monotone dynamical systems
\cite{halsmith,hirschsmith}, and use it to prove global convergence
in certain classes of CRNs where only very mild assumptions are made
on the kinetics.

The key geometrical insight on which the results are built was
provided in \cite{mierczynski} and generalised in
\cite{banajiangeli}. Stated very briefly, these results show that
sometimes the existence of an integral of motion in a strongly order-preserving dynamical system allows the construction of a Liapunov
function on each level set, which increases along nontrivial orbits.
This in turn has strong implications for the convergence of orbits.
While in general it may be an unusual conjunction of affairs to have
a first integral, strong monotonicity, and moreover integral and
order cone related in a particular way, it actually appears that this
situation is not uncommon in CRNs. However identifying when this
situation occurs is nontrivial, and explicit construction of families
of CRNs satisfying all these conditions becomes important.

The results at several points will be presented in considerably less
generality than possible in order to simplify the presentation and
highlight the key geometrical points, although where a theoretical
result allows greater generality, this may be mentioned.

\section{A convergence result}

The result presented in this section, and applied in subsequent
sections, is essentially derived from Theorem~2.4 in
\cite{banajiangeli}, with slight modification and specialisation for
our purposes. Note that Theorem~2.4 in \cite{banajiangeli} stated
that all orbits on a level set of the system in question converge to a unique
equilibrium, provided the equilibrium exists, while what was actually
proved was that all {\em bounded} orbits converge to the unique
equilibrium. As remarked in \cite{banajiangelierratum}, it remains an
open question whether the word ``bounded'' can be dropped from the
statement of the main result as erroneously done in
\cite{banajiangeli}. However for the purposes of this paper all
orbits will be bounded and only convergence of bounded orbits is
required. Note also that much of the difficulty in the proof of
Theorem~2.4 in \cite{banajiangeli} stemmed from the fact that the
integrals of motion considered were in general nonlinear whereas here
only the linear case is required.

Standard notions from convex geometry (as in
\cite{berman,nikaido,websterconvexity} for example) will be assumed.
Closed, convex and pointed cones in $\mathbb{R}^n$ define partial
orders on $\mathbb{R}^n$. Following \cite{berman}, cones in
$\mathbb{R}^n$ which are closed, convex, pointed and solid will be
referred to as {\bf proper}. Standard results in the theory of
monotone dynamical systems \cite{halsmith,hirschsmith} will also be
assumed.

{\bf Notation.} The symbols $<,>,\leq, \geq, \ll, \gg$ will refer to
the usual partial ordering of vectors in $\mathbb{R}^n$ derived from
the nonnegative orthant $\mathbb{R}^n_{\geq 0}$. So, given $a, b \in
\mathbb{R}^n$, $a \leq b$ means $b-a \in \mathbb{R}^n_{\geq 0}$, $a <
b$ means $b-a \in \mathbb{R}^n_{\geq 0} \backslash \{0\}$ and $a \ll
b$ means $b-a \in \mathrm{int}(\mathbb{R}^n_{\geq 0})$. When the
ordering is defined by some other closed, convex and pointed cone
$K$, the alternative symbols $\prec, \succ, \preceq, \succeq,
\llcurly, \ggcurly$ will be used. So, for example, $a \llcurly b$
means $b-a \in \mathrm{int}\,K$, and so forth.

Let $Y, K$ be proper cones in $\mathbb{R}^n$ with $K \supseteq Y$.
Define $K^{*}$ to be the dual cone to $K$, i.e., $K^* = \{y \in
\mathbb{R}^n\,|\, \langle y, k \rangle \geq 0\,\,\mbox{for all}\,\, k
\in K\}$. Consider a system
\begin{equation}
\label{eq1}
\dot x = F(x)
\end{equation}
on $Y$, and assume that (\ref{eq1}) defines a local semiflow $\phi$
on $Y$.

\subsection{Three conditions} We define three conditions on
(\ref{eq1}) and the associated semiflow which will be referred to as
Conditions~1,~2~and~3:
\begin{enumerate}
\item
$\phi$ is {\bf monotone} with respect to $K$, i.e., given $x, y
\in Y$ with $x \prec y$ and any $t > 0$ such that $\phi_t(x)$ and
$\phi_t(y)$ are defined, $\phi_t(x) \prec \phi_t(y)$. Moreover,
$\phi$ is {\bf strongly monotone} in $\mathrm{int}\,Y$ in the
following sense: given $x \prec y$ with at least one of $x$ or
$y$ in $\mathrm{int}\,Y$, then $\phi_t(x) \llcurly \phi_t(y)$ for
all $t > 0$ such that $\phi_t(x)$ and $\phi_t(y)$ are defined.
\item
The system has an {\bf increasing linear first integral}, namely
there exists a linear function $H \colon \mathbb{R}^n \to
\mathbb{R}$ such that (i) $\nabla H \in \mathrm{int}(K^*)$, and
(ii) for each $y \in Y$, $\langle \nabla H, F(y) \rangle = 0$.
Assume (w.l.o.g.) that $H(0) = 0$.
\item
$\phi$ has no limit points on $\partial\,Y \backslash \{0\}$.
\end{enumerate}

\subsection{Remarks on the conditions}

{\bf Condition 1.} Verifying monotonicity/ strong monotonicity of
$\phi$ is generally the most nontrivial step in applications. The
natural approach is via examination of the Jacobian matrix of the
system as will be spelled out in Section~\ref{secJac} below.

{\bf Condition 2.} For each $h \in \mathbb{R}$ define the level set
$S_h = \{y \in Y\,|\,H(y) = h\}$. The condition $\langle \nabla H,
F(y) \rangle = 0$ tells us that $H$ is constant along trajectories,
i.e., the sets $S_h$ are invariant. Choosing any set of basis vectors
for $(\nabla H)^\perp$ and writing these as the columns of a matrix
$\Gamma$ we can write $F(y) = \Gamma v(y)$ for some function $v \colon Y \to \mathbb{R}^{n-1}$. The implication of $\nabla H \in
\mathrm{int}(K^*)$ is as follows: given $x, y \in Y$ with $x\succ y$,
by convexity of $Y$, the line segment between $x$ and $y$ lies in
$Y$, and by integrating $\nabla H$ along this line segment we get
$H(x) > H(y)$. So $H$ is increasing with respect to the partial order
defined by $K$. This implies in turn that the level sets of $H$ are
unordered, and since $H(0) = 0$, $H(y) > 0$ for all $y \in
Y\backslash\{0\}$ (so $S_0 = \{0\}$). Since $H$ is linear and $Y$ is
a cone, clearly the range of $H$ is $[0,\infty)$, that is, $H \colon Y \to [0,\infty)$ is surjective.

{\bf Condition 3.} Condition $3$ is stronger than standard
``persistence'' assumptions which require only that limit sets of
{\em interior} points in $Y$ may not intersect the boundary. The
particular class of examples discussed here can be shown to fulfil
Condition 3; however making a weaker persistence assumption leads to
similar, slightly rephrased, convergence results to the ones
presented here. The choice to make a strong persistence assumption is
to reduce technical digressions and allow the main arguments to be
more transparent.

\subsection{The key theoretical result}

The basic theoretical result from which applications to CRNs will
follow is:

\begin{thm}
\label{mainthm}
Consider System (\ref{eq1}) with Conditions 1, 2 and 3 above. For
each $h \in [0, \infty)$, $S_h$ contains a unique equilibrium $e(h)$,
and $e(h)$ attracts all of $S_h$.
\end{thm}
\begin{pf}
Let $E = \{y \in Y\,|\, F(y) = 0\}$ be the equilibrium set of the
system. The main landmarks in the proof are presented below, but the
reader may refer to \cite{banajiangeli} to fill in the details.
\begin{enumerate}
\item
Any two equilibria $e_1, e_2$ must satisfy $e_1 \ggcurly e_2$ or
$e_1 \llcurly e_2$. It follows from Conditions~2~and~3 that $E
\cap
\partial Y = \{0\}$. Consequently given two distinct equilibria
$\bar{e}_1, \bar{e}_2$ at least one of these must lie in
$\interior{Y}$. As a result, if $\bar{e}_1, \bar{e}_2$ satisfy
$\bar{e}_1 \prec \bar{e}_2$, then by strong monotonicity in
$\interior{Y}$, $\bar{e}_1 \llcurly \bar{e}_2$. Now suppose there
are equilibria $\bar{e}_1, \bar{e}_2$ which are unordered by
$\prec$.  By repeating the geometrical arguments as in the proof
of Lemma~5.11 of \cite{banajiangeli} we obtain the existence of
an equilibrium $\tilde{e}$ such that $\tilde{e} \prec \bar{e}_1$
(say) but not $\tilde{e} \llcurly \bar{e}_1$, contradicting the
fact that $\tilde{e} \prec \bar{e}_1$ implies $\tilde{e} \llcurly
\bar{e}_1$.
\item
Each $S_h$ is a compact, convex set: by linearity of $H$, $S_h$
is the intersection of a hyperplane with $Y$ and so convex, while
boundedness follows from the condition that $\nabla H \in
\mathrm{int}(K^*)$ (see Lemma~3.9 in \cite{banajiangeli}). Thus
each level set contains an equilibrium. Uniqueness of the
equilibrium follows from the fact that $S_h$ is unordered,
whereas from 1, two equilibria $e_1$ and $e_2$ must satisfy $e_1
\ggcurly e_2$ or $e_1 \llcurly e_2$.
\item
The map $e \colon [0, \infty) \to E$ which associates to each $h$
the unique equilibrium $e(h)$ on $S_h$ is an order-preserving
homeomorphism. That $e$ is an order-preserving bijection with
continuous inverse (namely $\left. H\right|_E$) is immediate from
parts 1 and 2, and linearity of $H$. That $e$ is continuous
follows essentially from the fact that $E$ is closed (see
Lemma~5.12 in \cite{banajiangeli}).
\item
There is a continuous scalar function $L \colon Y \to [0,
\infty)$, increasing along nontrivial orbits in
$\mathrm{int}\,Y$, and which attains a global maximum on $S_h$ at
$e(h)$.  The proofs of these facts follow the related results in
\cite{banajiangeli} (Lemmas~5.13 to 5.18).  The value of $L(y)$
is the value of $H$ at the unique point $Q(y)$ where $y -
\partial K$ intersects $E$. Lemmas~5.13 to 5.17 in
\cite{banajiangeli} carry over to the present case immediately.
We claim that $L$ is increasing along orbits of those $y \notin
E$ for which $L(y) > 0$.  Indeed, consider any such $y$.  By the
definition of $Q(y)$, $Q(y) \prec y$.  Since $L(y) = H(Q(y)) >
0$, we have $Q(y) \ne 0$, so, by Condition 3, $Q(y) \in
\interior{Y}$.  Strong monotonicity in $\interior{Y}$ implies
that for any $t > 0$, $\phi_{t}(Q(y)) = Q(y) \llcurly
\phi_{t}(y)$. Consider any $e \in E$.  If $e \preceq Q(y)$ then
$e \llcurly \phi_{t}(y)$, i.e., $e \ne Q(\phi_{t}(y))$.  So
$Q(\phi_{t}(y)) \succ Q(y)$, and thus $L(\phi_{t}(y)) > L(y)$,
which proves the claim.  It remains to notice that $L(y) > 0$ for
any $y \in \interior{Y}$, since otherwise $y \in
\partial K$, which contradicts $\interior{Y} \subseteq
\interior{K}$. We remark that $L$ is well defined at points on
$\partial Y$, but may not be (strictly) increasing if $L(y) = 0$
(this is possible if $y \in
\partial Y \cap \partial K$) as the line segment $[0, y]$ lies in
$\partial Y$ where we do not necessarily have strong
monotonicity.
\item
If $h = 0$ it is trivial that $e(h)$ attracts all of $S(h)$ as
$S_h = e(h) = \{0\}$. So assume $h \not = 0$. Since $S_h$ is
compact, each $x \in S_h$ has a nonempty $\omega$-limit set,
$\omega(x)$. The Liapunov function in $\mathrm{int}\,Y$, strictly
increasing along nontrivial trajectories and achieving a maximum
at $e(h)$, guarantees that $\omega(x)$ does not intersect
$\mathrm{int}\,Y\backslash\{e(h)\}$. So $\omega(x) \subseteq \partial Y \cup \{e(h)\}$. But by assumption $\omega(x)$ does
not intersect $\partial\,Y\backslash\{0\}$ and $0 \not \in S_h$,
so $\omega(x) = \{e(h)\}$.
\end{enumerate}

\qquad \qed

\end{pf}

{\bf Remark.} For $h \not = 0$ the equilibrium $e(h)$ defined in the
above proof clearly lies in $\mathrm{int}\,Y$ and so in the relative
interior of $S_h$.

\subsection{Jacobian conditions for monotonicity and strong
monotonicity}
\label{secJac}

In general, the greatest difficulty in applying Theorem~\ref{mainthm}
is verifying monotonicity of a local semiflow in $Y$ and strong
monotonicity in $\mathrm{int}\,Y$ in the sense described in Condition
1. For this reason it is useful to have conditions which can more
easily be checked for a system.

{\bf Terminology.} Let $J$ be a square matrix. Given a closed, convex
and pointed cone $K \subseteq \mathbb{R}^n$, we will say that $J$ is
{\bf $K$-positive} if $y \in K$ implies $Jy \in K$. A $K$-positive
matrix $J$ is {\bf $K$-irreducible} if the only faces $F$ of $K$ for
which $J(F) \subseteq F$ are $\{0\}$ and $K$. $J$ is {\bf
$K$-quasipositive} if given any $y \in K$, there exists $\alpha \in
\mathbb{R}$ such that $Jy + \alpha y \in K$. Suppose that there
exists $\alpha \in \mathbb{R}$ such that $J + \alpha I$ is
$K$-positive and $K$-irreducible, then we say that $J$ is {\bf
strictly $K$-quasipositive}.

The following is a slight adaptation of results in \cite{hirschsmith}
for our purposes.

\begin{lem}
\label{hslem}
Consider a proper cone $K \subseteq \mathbb{R}^n$, some open set $U
\subseteq \mathbb{R}^n$ , and a $C^1$ vector field $f \colon U \to
\mathbb{R}^n$ with Jacobian matrix $Df$. Let $X \subseteq U$ be some
convex domain, positively invariant under the local flow $\phi_U$
defined by $f$, and let $\phi$ be the induced local semiflow on $X$.
Assume that $Df$ is $K$-quasipositive in $X$. Consider some $x_0, x_1
\in X$ with $x_0 \prec x_1$. Then $\phi_t(x_0) \prec \phi_t(x_1)$ for
each $t > 0$ such that $\phi_t(x_0), \phi_t(x_1)$ exist. If there
exists $y_0$ on the line segment joining $x_0$ to $x_1$ such that
$Df(y_0)$ is strictly $K$-quasipositive, then $\phi_t(x_0) \llcurly
\phi_t(x_1)$ for each $t > 0$ such that $\phi_t(x_0), \phi_t(x_1)$
exist.
\end{lem}
\begin{pf}
By basic results on monotone dynamical systems (see Section~3.1 in
\cite{hirschsmith}), $K$-quasipositivity implies monotonicity: given
$x_0, x_1 \in X$ with $x_0 \prec x_1$, $\phi_t(x_0) \prec
\phi_t(x_1)$ for each $t > 0$ such that $\phi_t(x_0), \phi_t(x_1)$
exist. Note that if $Df(y_0)$ is strictly $K$-quasipositive, then for
each $x \in \partial K\backslash\{0\}$, there exists $\nu \in K^{*}$
such that $\nu(x) = 0$ and $\nu(Df(y_0)x) > 0$ (see Proposition 3.10
of \cite{hirschsmith}). The final statement now follows from Lemma
3.7 in \cite{hirschsmith}. \qquad \qed
\end{pf}

{\bf Remark.} Since we are discussing the Jacobian matrix it reduces
technicalities to assume that the vector field is defined on some
open set containing the positively invariant set $X$ in which we are
interested.

\section{Some background material}

From now on we fix the state space of interest as
$\mathbb{R}^{m}_{\geq 0}$, the nonnegative orthant in
$\mathbb{R}^{m}$ and focus on systems with a linear first integral.
Let $\Gamma$ be an $m \times n$ matrix, and $v \colon U \to
\mathbb{R}^{n}$ be a $C^1$ function where $U$ is some open
neighbourhood of $\mathbb{R}^{m}_{\geq 0}$. We consider the system
\begin{equation}
\label{maineqn}
\dot x = \Gamma v(x)
\end{equation}
defining a local flow $\tilde \phi$ on $U$. We assume that
$\mathbb{R}^{m}_{\geq 0}$ is (positively) invariant under $\tilde
\phi$ and so get a local semiflow $\phi$ on $\mathbb{R}^{m}_{\geq
0}$.

Dynamical systems arising from chemical reaction networks typically
take the form (\ref{maineqn}). In this context, the matrix $\Gamma$
is termed the \textbf{stoichiometric matrix} of the system, and the
intersection of any coset of $\mathrm{Im}\,\Gamma$ with
$\mathbb{R}^{m}_{\geq 0}$ is called a \textbf{stoichiometry class}.
Clearly stoichiometry classes are positively invariant under $\phi$.
The function $v$ describes the rates of reaction. A variety of common
forms for $v$ exist in the experimental and modelling literature, all
of which satisfy the weak assumptions we present later
(Section~\ref{seckinetics}).

\subsection{Qualitative classes and sets of matrices}

A real matrix $M$ determines the following three sets of matrices:
\begin{enumerate}
\item
The qualitative class $\mathcal{Q}(M)$ \cite{brualdi} of all
matrices with the same sign pattern as $M$. Explicitly, $N \in
\mathcal{Q}(M)$ if $N$ has the same dimensions as $M$ and
satisfies $\mathrm{sign}(N_{ij}) = \mathrm{sign}(M_{ij})$, i.e.,
$M_{ij} > 0 \Rightarrow N_{ij} > 0$, $M_{ij} < 0 \Rightarrow
N_{ij} < 0$ and $M_{ij} = 0 \Rightarrow N_{ij} = 0$.
\item
The set $\mathcal{Q}_0(M) \equiv \mathrm{cl}(\mathcal{Q}(M))$,
the closure of $\mathcal{Q}(M)$. $\mathcal{Q}_0(M)$ consists of
all matrices $N$ such that $M_{ij} > 0 \Rightarrow N_{ij} \geq
0$, $M_{ij} < 0 \Rightarrow N_{ij} \leq 0$ and $M_{ij} = 0
\Rightarrow N_{ij} = 0$.
\item
The set $\mathcal{Q}_1(M)$ consisting of all matrices $N$ such
that $M_{ij} > 0 \Rightarrow N_{ij} \geq 0$, $M_{ij} < 0
\Rightarrow N_{ij} \leq 0$.
\end{enumerate}
The following fact will often be applied without comment: if $c \in
\mathcal{Q}_0(b)$ and $b \in \mathcal{Q}_1(a)$, then $c \in
\mathcal{Q}_1(a)$. The next lemma can easily be checked directly.
\begin{lem}
\label{Q1lem}
Let $a, b \in \mathbb{R}^n$, and let $D$ be an $n \times n$ diagonal
matrix with positive diagonal entries. Suppose $D a, D b$ are
$(-1,0,1)$-vectors. Then $a+b \in \mathcal{Q}_1(a) \cap
\mathcal{Q}_1(b)$. \qquad \qed
\end{lem}

\subsection{Associating digraphs with CRNs}

A DSR graph \cite{banajicraciun2} is a signed, labelled, bipartite,
multidigraph, examination of which can lead to a number of
conclusions about the behaviour of a system of differential
equations. In this paper we will be interested only in whether a
given digraph is strongly connected, namely whether there is a
(directed) path from each vertex in the digraph to each other vertex.
So, following \cite{banajiburbanks}, we present a closely related but
simplified construction. Given an $m \times n$ matrix $A$ and an $n
\times m$ matrix $B$, we construct a bipartite digraph $G_{A, B}$ on
$m + n$ vertices as follows: associate a set of $m$ vertices $u_1,
\ldots, u_m$ with the rows of $A$; associate another $n$ vertices
$v_1, \ldots, v_n$ with the columns of $A$; add the arc $u_iv_j$ iff
$A_{ij} \not = 0$; add the arc $v_ju_i$ iff $B_{ji} \not = 0$.

Such a digraph $G(x)$ can be associated with System~(\ref{maineqn})
at each point in $x \in \mathbb{R}^{m}_{\geq 0}$ as follows: define
$V(x) \equiv \left[\partial v_i(x)/ \partial x_j\right]$, and let
$G(x) \equiv G_{\Gamma, V(x)}$. \\

{\bf Remark.} It is in general possible to factorise a given function
in different ways, say $F(x) = \Gamma v(x) = \tilde{\Gamma}
\tilde{v}(x)$. These different factorisations result in different
digraphs. Thus the digraph associated with System~(\ref{maineqn}) is
not strictly associated with the vector field itself, but with a
particular factorisation -- see \cite{banajicraciun2} for more
discussion of this. In the case of CRNs, there is a natural
factorisation, where $\Gamma$ is the stoichiometric matrix.

\subsection{The $n$-cube}

Given a matrix $M$, $M_k$ will refer to the $k$th column of $M$, and $M^k$ to the $k$th row of $M$.\\

\begin{figure}[h]
\begin{center}
\begin{tikzpicture}[scale=2]
	 \tikzstyle{vertex}=[circle,minimum size=12pt,inner sep=0pt]
	 \tikzstyle{selected vertex} = [vertex, fill=red!24]
	 \tikzstyle{selected edge} = [draw,line width=5pt,-,red!50]
	 \tikzstyle{edge} = [draw,line width=0.04cm,,-,black]
	 \node[vertex] (v0) at (0,0) {$1$};
	 \node[vertex] (v1) at (0,1) {$5$};
	 \node[vertex] (v2) at (1,0) {$2$};
	 \node[vertex] (v3) at (1,1) {$6$};
	 \node[vertex] (v4) at (0.23, 0.4) {$3$};
	 \node[vertex] (v5) at (0.23,1.4) {$7$};
	 \node[vertex] (v6) at (1.23,0.4) {$4$};
	 \node[vertex] (v7) at (1.23,1.4) {$8$};
	 \draw[->, thick, color=black!60] (v2) -- (1.5,0);\node at (1.63,-0.01) {$x_1$};
	 \draw[->, thick, color=black!60] (v4) -- (0.39,0.7);\node at (0.45,0.8) {$x_2$};
	 \draw[->, thick, color=black!60] (v1) -- (0,1.5);\node at (0,1.6) {$x_3$};
	 \draw[edge] (v0) -- (v4) -- (v5) -- (v1) -- (v0);
	 \draw[edge] (v2) -- (v6) -- (v7) -- (v3) -- (v2);
	 \draw[edge] (v4) -- (v6) -- (v7) -- (v5) -- (v4);
	 \draw[-, line width=0.1cm, color=white] (v0) -- (v1) -- (v3) -- (v2) -- (v0);

	 \draw[edge] (v0) -- (v1) -- (v3) -- (v2) -- (v0);

 \end{tikzpicture}
\end{center}
\caption{\label{cubefig} The $3$-cube in $\mathbb{R}^3$ with vertex
$n$ corresponding to column $n$ of $B$. $B_{kn}$ is the $(4-k)$th entry in the binary representation of $n-1$. Observe that if $i-1 <
2^{k-1} \mod 2^k$, then $B_{i+2^{k-1}}-B_i = e_k$ corresponding to
the geometrical fact that each basis vector is associated with $4$ parallel edges of the cube defined by $B$. For example $B_{3+2^2} - B_3 = B_7 - B_3 =
(0,1,1)^T - (0,1,0)^T = (0, 0, 1)^T = e_3$. }
\end{figure}
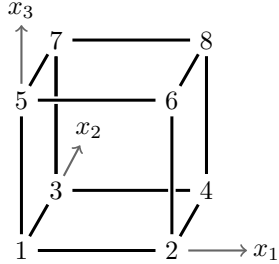

Fix $n \geq 2$ and define the $n \times 2^{n}$ matrix $B$ as the
$(0,1)$-matrix whose $j$th column read vertically upwards as a binary
string gives the number $j-1$ (so $\sum_{i=1}^n2^{i-1}B_{ij} = j-1$).
For example, if $n = 3$:
\[
B =
\left(\begin{array}{cccccccc}
0&1&0&1&0&1&0&1
\\
0&0&1&1&0&0&1&1
\\
0&0&0&0&1&1&1&1\end{array}\right).
\]
Note that the columns of $B$ are just the vertices of the
$n$-dimensional hypercube, and two vertices are adjacent if they
differ in exactly one coordinate. The choice of ordering is
important: given $i < j$, $B_i$ and $B_j$ are adjacent iff (i) there
exists $ k \geq 0$ such that $j = i + 2^{k}$ and (ii) $i-1 < 2^{k}
\mod 2^{k+1}$. The following is immediate:
\begin{lem}
\label{lemB}
Given $i \in \{1, \ldots, 2^n\}$ and $k \in \{1, \ldots, n\}$, if
$i-1 \geq 2^{k-1} \mod 2^k$, then $B_i - B_{i-2^{k-1}} = e_k$ and
otherwise $B_{i+2^{k-1}}-B_i = e_k$. \qquad \qed
\end{lem}
The geometrical meaning of Lemma~\ref{lemB} is illustrated for $n =
3$ in Figure~\ref{cubefig}.

\subsection{Cubic cones}

We now construct a class of cones which we refer to as ``cubic
cones''. First we show that no binary vector is in the convex hull of
other, different, binary vectors. Geometrically this is obvious:
these binary vectors are precisely the vertices of the $n$-cube.
\begin{lem}
\label{binconvex}
Let $B$ be any $(0, 1)$-matrix without repeated columns and let $p$
be a nonnegative vector with $\sum_ip_i = 1$. If $Bp = q$ and $q$ is
a $(0, 1)$-vector then $p = e_j$ for some $j$.
\end{lem}
\begin{pf}
We look for solutions to $q \equiv Bp$ with $q$ a binary vector and
$\sum_i p_i = 1$. For each $i$,
let $S_i^0$ be the set of indices such that $k \in S_i^0
\Leftrightarrow B_{ik} = 0$. Similarly define $S_i^1$ via $k \in
S_i^1 \Leftrightarrow B_{ik} = 1$. Then:
\begin{enumerate}
\item
If $q_i = \sum_kB_{ik}p_k = 1$, then $\sum_{k \in S_i^1}p_k = 1$,
so $\sum_{k \in S_i^0}p_k = 0$, so $p_k = 0$ for each $k \in
S_i^0$.
\item
If $q_i = \sum_kB_{ik}p_k = 0$, then $\sum_{k \in S_i^1}p_k = 0$,
so $p_k = 0$ for each $k \in S_i^1$.
\end{enumerate}
Since $p$ is not the zero vector, there exists $j$ such that $p_j
\not = 0$. Then for each $i$ such that $q_i = 1$, $j \not \in S_i^0$,
i.e., $j \in S_i^1$, i.e., $B_{ij}=1$; similarly for each $i$ such
that $q_i = 0$, $j \not \in S_i^1$, i.e., $j \in S_i^0$, i.e.,
$B_{ij} = 0$. We have $q_i = B_{ij}$, i.e., $q$ is just the $j$th
column of $B$. If there exists $l \not = j$ such that $p_l \not = 0$,
then we get $q_i = B_{il}$. Since $B$ has no repeated columns this is
impossible, so $p_l = 0$ for $l \not = j$. Since $\sum_ip_i = 1$, $p
= e_j$. \qquad \qed
\end{pf}

The construction in the next lemma is crucial to what is to follow:
\begin{lem}
\label{cubicconeslem}
Let $m > n$ and consider an $m \times n$ matrix $\Gamma$ of rank $n$.
Let $c \in \mathbb{R}^{m}$ be some vector not in
$\mathrm{Im}\,\Gamma$. Define $\mathbf{1}$ to be the vector in
$\mathbb{R}^{2^{n}}$ each of whose entries is $1$ and define $\Lambda
\equiv c\mathbf{1}^\mathrm{T} + \Gamma B$. With $\Lambda$ defined in
this way,
\begin{enumerate}
\item[{\rm (i)}]
$K \equiv \{\Lambda z\,|\, z \geq 0\}$ is a closed, convex and
pointed cone in $\mathbb{R}^{m}$;
\item[{\rm (ii)}]
The set of those $r \in \mathrm{ker}\,\Gamma^\mathrm{T}$ for
which $r^\mathrm{T} c > 0$ is nonempty and is contained in
$\mathrm{int}(K^*)$;
\item[{\rm (iii)}]
The nonnegative multiples of the $2^n$ columns of $\Lambda$ are
precisely the closed one-dimensional faces of $K$.
\end{enumerate}
\end{lem}

\begin{pf}
That $K$ is a closed, convex cone is immediate from the definition.
Since $c \not \in \mathrm{Im}\,\Gamma$, there exists $r \in
\mathrm{ker}\,\Gamma^\mathrm{T}$ such that $r^\mathrm{T} c > 0$. For
any $p > 0$ and any such $r$, $r^\mathrm{T}\Lambda p =
r^\mathrm{T}c\mathbf{1}^\mathrm{T}p = r^\mathrm{T}c\sum p_i > 0$. So
$r \in \mathrm{int}(K^*)$. Since $\mathrm{int}(K^*)$ is nonempty, $K$
is pointed.

The third statement fails to be true only if some column of $\Lambda$
can be written as a nonnegative combination of other columns of
$\Lambda$. Suppose $\Lambda e_k = \Lambda p$ where $0 \leq p$. Then:
\[
(c\mathbf{1}^\mathrm{T} + \Gamma B)e_k = (c\mathbf{1}^\mathrm{T} +
\Gamma B)p
\]
i.e.,
\[
c(1 - \sum p_i) = \Gamma B (p-e_k)
\]
Left-multiplying both sides by $r^T$ gives $r^Tc(1 - \sum p_i) = 0$,
so $\sum p_i = 1$ and $\Gamma B (p-e_k) = 0$. Since  $\Gamma$ has
rank $n$, $\mathrm{ker}\,\Gamma = \{0\}$, and so $B(p-e_k) = 0$. By
Lemma~\ref{binconvex}, $p=e_k$. So each column of $\Lambda$ is an
extremal vector of $K$. \qquad \qed
\end{pf}

{\bf Geometrical interpretation.} Since $\Gamma$ defines a linear map
from $\mathbb{R}^n$ to $\mathbb{R}^m$ with trivial kernel and $B$ has
rank $n$, the columns of $\Gamma B$ are the extremal points on an
$n$-dimensional parallelotope in $\mathbb{R}^m$.
%The adjacency relationships are preserved by a linear map.
Adding on the constant vector $c$ simply translates these points, so
the columns of $\Lambda$ also define an $n$-parallelotope in
$\mathbb{R}^m$. Ensuring that the constant vector $c$ does not lie in
$\mathrm{Im}\,\Gamma$ means that the cone $K(\Lambda) \equiv
\{\Lambda z\,|\, z \geq 0\}$ is pointed, and the convex hull of
columns of $\Lambda$, $\mathrm{conv}(\Lambda)$, is a cross-section of
this cone. As $\mathrm{conv}(\Lambda)$ is simply the $n$-cube after
an affine mapping we call these cones {\bf cubic cones} and the
associated partial orders {\bf cubic orders}.

Obviously, the incidence structure of $\mathrm{conv}(\Lambda)$ is
precisely that of the cube, $\mathrm{conv}(B)$, namely the convex
hull of any $2^k$ vectors $\Lambda_{r_1}, \ldots, \Lambda_{r_{2^k}}$
is a $k$-dimensional face of $\mathrm{conv}(\Lambda)$ if and only if
the convex hull of the corresponding vectors $B_{r_1}, \ldots,
B_{r_{2^k}}$ is a $k$-dimensional face of the $n$-cube.

Each vector $\Gamma_k$ is associated with $2^{n-1}$ edges of $\mathrm{conv}(\Lambda)$ corresponding to $2^{n-1}$ parallel edges of $\mathrm{conv}(B)$. More
precisely, by Lemma~\ref{lemB}, if $i-1 \geq 2^{k-1} \mod 2^k$ then
\[
\Lambda_i - \Lambda_{i-2^{k-1}} = c + \Gamma B_i - (c + \Gamma
B_{i-2^{k-1}}) = \Gamma (B_i - B_{i-2^{k-1}}) = \Gamma_k
\]
while if $i-1 < 2^{k-1} \mod 2^k$ then
\[
\Lambda_{i+2^{k-1}} - \Lambda_i = c + \Gamma B_i - (c + \Gamma
B_{i-2^{k-1}}) = \Gamma (B_{i+2^{k-1}} - B_i) = \Gamma_k
\]
So for any $k \in \{1, \ldots, n\}$ and any $i \in \{1, \ldots,
2^n\}$ we define
\[
j(i, k) \equiv \left\{\begin{array}{ll}i-2^{k-1} & \mbox{(if $i-1
\geq 2^{k-1} \mod 2^k$)}\\i+2^{k-1}&
\mbox{(otherwise)}\end{array}\right.
\]
and can always write for any $i, k$
\[
\Gamma_k = \mathrm{sign}(j(i, k) - i)(\Lambda_{j(i,k)} -\Lambda_i).
\]

\section{Systems preserving cubic orders}

We proceed to construct a class of systems to which
Theorem~\ref{mainthm} can be applied. We do this by prescribing {\it
a priori} the geometric nature of the order to be preserved and then
constructing systems which preserve such orders. Via the next two
theorems we provide a set of easily checkable conditions which ensure
that System~(\ref{maineqn}) has the behaviour described in
Theorem~\ref{mainthm}. Later we will provide examples of CRNs which
fulfil these conditions.

\begin{thm}
\label{mainthm0}
Consider System~(\ref{maineqn}) with $\Gamma$ now an $(n+1) \times n$
matrix of rank $n$. Assume:
\begin{enumerate}
\item
For some $c \not = 0$, $\Lambda \equiv c\mathbf{1}^\mathrm{T} +
\Gamma B$ satisfies:
\begin{enumerate}
\item
There exists a nonnegative right-inverse to $\Lambda$, i.e.,
a nonnegative matrix $P$ such that $\Lambda P = I$.
\item
There exists a diagonal matrix $D$ with positive diagonal
entries such that $D\Lambda$ is a $(-1,0,1)$-matrix
\end{enumerate}
\item
For each $x \in \mathbb{R}^{n+1}_{\geq 0}$, $V(x) \equiv
\left[\partial v_i(x)/ \partial x_j\right] \in
\mathcal{Q}_0(-\Gamma^\mathrm{T})$
\item
For each $x \in \mathrm{int}(\mathbb{R}^{n+1}_{\geq 0})$, the
digraph $G_{\Gamma, V(x)}$,  is strongly connected.
\end{enumerate}
Define
\[
K = K(\Lambda) \equiv \{\Lambda z\,|\, z \geq 0\}.
\]
Then:

\begin{enumerate}
\item
$\mathbb{R}^{n+1}_{\geq 0} \subseteq K$.
\item
$c \not \in \mathrm{Im}\,\Gamma$. $\mathrm{int}(K^*)$ is nonempty
and $K$ is pointed.
\item
Each column of $\Lambda$ is an extremal vector of $K$.
\item
For each $x \in \mathbb{R}^{n+1}_{\geq 0}$, $\Gamma V(x)$ is
$K$-quasipositive, while for each $x \in
\mathrm{int}(\mathbb{R}^{n+1}_{\geq 0})$, $\Gamma V(x)$ is
strictly $K$-quasipositive.
\end{enumerate}
\end{thm}

Before proving Theorem~\ref{mainthm0} we show how when it is combined
with Theorem~\ref{mainthm} we get a global convergence result.

\begin{thm}
\label{applicationthm}
Consider System~(\ref{maineqn}) with the assumptions in
Theorem~\ref{mainthm0}. Assume further that
$\partial\,\mathbb{R}^{n+1}_{\geq 0}\backslash\{0\}$ contains no
limit points of the local semiflow $\phi$.  Then there exists a
linear function $H \colon \mathbb{R}^{n+1}_{\geq 0}\to \mathbb{R}$
whose level sets \[ S_h = \{y \in \mathbb{R}^{n+1}_{\geq 0}\,|\, H(y)
= h\}.
\]
are invariant under $\phi$. Moreover if for some $h$, $S_h$ is
nonempty, it contains a unique equilibrium $e(h)$ which attracts all
of $S_h$.
\end{thm}
\begin{pf}
The cone $K$ in Theorem~\ref{mainthm0} is clearly closed and convex.
That it is pointed and solid follow from results 1 and 2 of the
theorem. So $K$ is a proper cone. Result 4 of Theorem~\ref{mainthm0}
gives that the Jacobian matrix of the system is $K$-quasipositive
everywhere in $\mathbb{R}^{n+1}_{\geq 0}$ and strictly
$K$-quasipositive in $\mathrm{int}(\mathbb{R}^{n+1}_{\geq 0})$. By
Lemma~\ref{hslem}, $K$-quasipositivity implies that given $x_0 \prec
x_1$, $\phi_t(x_0) \prec \phi_t(x_1)$ (for each $t > 0$ such that
$\phi_t(x_0), \phi_t(x_1)$ exist); further, strict
$K$-quasipositivity in $\mathrm{int}(\mathbb{R}^{n+1}_{\geq 0})$
implies that if either of $x_0$ or $x_1$ lies in
$\mathrm{int}(\mathbb{R}^{n+1}_{\geq 0})$ then $\phi_t(x_0) \llcurly
\phi_t(x_1)$. Thus Condition~1 of Theorem~\ref{mainthm} is satisfied.

Since by Theorem~\ref{mainthm0}(2), $c \not \in \mathrm{Im}\,\Gamma$,
Lemma~\ref{cubicconeslem}(ii) allows us to choose $r \in
\mathrm{ker}\,\Gamma^\mathrm{T}$ such that $r^\mathrm{T} c > 0$.
Define $H(x) = r^\mathrm{T} x$. $H$ is by definition linear, with $H(0) = 0$, and clearly
\[
\langle \nabla H, \Gamma v(y) \rangle = r^T\Gamma v(y) = 0
\]
so that level sets of $H$ are invariant.
Lemma~\ref{cubicconeslem}(ii) provides further that $\nabla H = r \in
\mathrm{int}(K^*)$, so the integral is increasing. Thus Condition~2
is satisfied.

Finally, Condition 3 is assumed. Thus the assumptions of
Theorem~\ref{mainthm} are fulfilled and so Theorem~\ref{mainthm} can
be applied. \qquad \qed
\end{pf}

A somewhat more general statement than the following is proved as
Theorem~1 in \cite{banajiburbanks}.

\begin{lem}
\label{lemDSR}
Let $K \subseteq \mathbb{R}^n$ be a closed, convex and pointed cone,
$A$ an $n \times m$ matrix, and $B'$ an $m \times n$ matrix. Suppose
that $\mathrm{Im}\,A \not \subseteq \mathrm{span}\,F$ for any
nontrivial face $F$ of $K$, and that $AB$ is $K$-quasipositive for
each $B \in \mathcal{Q}_0(B')$. Then whenever $G_{A, B}$ is strongly
connected, $AB$ is strictly $K$-quasipositive.
\end{lem}

\subsection{Proof of theorem~\ref{mainthm0}}

The lemmas presented so far make the proof now straightforward.

\begin{pot2}

\textbf{1. $\mathbb{R}^{n+1}_{\geq 0} \subseteq K$.} Since $\Lambda
Pe_i = e_i$ and $Pe_i$ is nonnegative, $e_i \in K$ for each $i$ and
so by convexity $\mathbb{R}^{n+1}_{\geq 0} \subseteq K$.

\textbf{2. $c \not \in \mathrm{Im}\,\Gamma$. There exists $r \in
\interior(K^*)$ and $K$ is pointed.} From the previous part,
$\mathbb{R}^{n+1} = K-K$. One consequence is that, regarded as a map
from $\mathbb{R}^{2^n}$ to $\mathbb{R}^{n+1}$, $\Lambda$ is
surjective. It follows that $c \not \in \mathrm{Im}\,\Gamma$:
otherwise, since all vectors in $\mathrm{Im}\,\Lambda$ are of the
form $\lambda c + \Gamma b$ (for some scalar $\lambda$ and some
vector $b$), each vector in $\mathrm{Im}\,\Lambda$ could be written
$\Gamma d$ for some vector $d$; however $\Gamma$ only has rank $n$
contradicting the fact that $\Lambda$ is surjective.

Since $c \not \in \mathrm{Im}\,\Gamma$, by
Lemma~\ref{cubicconeslem}(ii) there exists $r \in
\mathrm{ker}\,\Gamma^\mathrm{T}$ for which $r^\mathrm{T} c > 0$, and
$r \in \mathrm{int}(K^*)$.  Thus $K$ is pointed. Note that since
$\mathbb{R}^{n+1}_{\geq 0} \subseteq K$, it follows that
$\interior(K^*) \subseteq \interior(\mathbb{R}^{n+1}_{\geq 0})$, so
$r \gg 0$.

{\bf 3. Each column of $\Lambda$ is an extremal vector of $K$.} This
follows from Lemma~\ref{cubicconeslem}(iii).

{\bf 4a. For each $y \in \mathbb{R}^{n+1}_{\geq 0}$, $\Gamma V(y)$ is
$K$-quasipositive.} Fix $y$ and let $V = V(y)$. Recall that for $i
\in \{1, \ldots, 2^n\}$ and $k \in \{1, \ldots, n\}$ there exists
$j(i, k) \in \{1, \ldots, 2^n\}$ such that
\begin{equation*}
\Gamma_k = \mathrm{sign}(j(i, k) - i)(\Lambda_{j(i,k)} -\Lambda_i).
\end{equation*}
Define $\alpha_{ik} = |V^k\Lambda_i|$,  $\alpha_i = \sum_k
\alpha_{ik}$ and $\alpha = \max_i\,\alpha_i$. We have
\[
\Gamma V \Lambda_i = \sum_k \Gamma_k V^k \Lambda_i.
\]
If $i > j(i, k)$, then $\Gamma_k = \Lambda_{i}-\Lambda_{j(i,k)}$, and
by Lemma~\ref{Q1lem} we have $\Gamma_k \in
\mathcal{Q}_1(\Lambda_{i})$. Since $-V^k \in
\mathcal{Q}_0(\Gamma_k^T)$, $-V^k \in \mathcal{Q}_1(\Lambda_i^T)$, so
$-V^k\Lambda_i \geq 0$. Similarly, $i < j(i, k)$ implies
$V^k\Lambda_i \geq 0$. So, in either case
\begin{eqnarray*}
\Gamma_k V^k \Lambda_i + \alpha_{ik} \Lambda_i & = &
\mathrm{sign}(j(i, k) - i)(\Lambda_{j(i,k)} -\Lambda_i)V^k \Lambda_i
+ \alpha_{ik} \Lambda_i
\\
& = & (\Lambda_{j(i,k)} -\Lambda_i)\alpha_{ik} + \alpha_{ik}
\Lambda_i  = \alpha_{ik}\Lambda_{j(i,k)}.
\end{eqnarray*}
So
\begin{equation}
\label{eqGV0}
\Gamma V \Lambda_i + \alpha_i \Lambda_i = \sum_k (\Gamma_k V^k
\Lambda_i + \alpha_{ik} \Lambda_i) = \sum_k
\alpha_{ik}\Lambda_{j(i,k)} \in K
\end{equation}
since $\alpha_{ik} \geq 0$ for each $i$. Given any $x = \sum_i t_i
\Lambda_i\in K$ ($t_i \geq 0$ for each $i$):
\begin{eqnarray}
\Gamma V x + \alpha x & = & \sum_i t_i \left(\alpha \Lambda_i +
\Gamma V \Lambda_i\right) \nonumber
\\
\label{eqGV}
& = & \sum_i t_i \left((\alpha - \alpha_i) \Lambda_i + \Gamma V
\Lambda_i + \alpha_i \Lambda_i \right)
\end{eqnarray}
Since $\alpha - \alpha_i \geq 0$ for each $i$ and $\Gamma V \Lambda_i
+ \alpha_i \Lambda_i \in K$, $\Gamma V x + \alpha x \in K$.

{\bf 4b. For each $y \in \mathrm{int}(\mathbb{R}^n_{\geq 0})$,
$\Gamma V(y)$ is strictly $K$-quasipositive.} By 4(a) we already know
that $\Gamma V$ is $K$-quasipositive for all $V$ in
$\mathcal{Q}_0(-\Gamma^T)$. By assumption, for each $y \in
\mathrm{int}(\mathbb{R}^{n+1}_{\geq 0})$, the digraph $G_{\Gamma,
V(y)}$,  is strongly connected. In order to apply Lemma~\ref{lemDSR}
we need to confirm that $\mathrm{Im}\,\Gamma \not \subseteq
\mathrm{span}\,F$ for any nontrivial face $F$ of $K$.

Consider a nontrivial face $F$ of $K$. We put the columns $\Lambda_i
\in F$ into a matrix $\Lambda_F$. Suppose $\mathrm{Im}\,\Gamma
\subseteq \mathrm{Im}\,\Lambda_F$, i.e., there exists a matrix $T$
such that $\Gamma = \Lambda_FT$. Since $\Lambda_F$ simply consists of
a subset of the columns of $\Lambda$, for any $i$ there exists $j(i)$
such that $\Lambda_Fe_i = c + \Gamma B_{j(i)}$, and as $\Gamma =
\Lambda_FT$, $c = \Lambda_F(e_i - T B_{j(i)})$. Thus $\Lambda =
\Lambda_F[(e_i - T B_{j(i)})\mathbf{1}^T + T B]$. This shows that
$\mathrm{Im}\,\Lambda \subseteq \mathrm{Im}\,\Lambda_F$,
contradicting the fact that $F$ is a nontrivial face of $K$. Thus
$\mathrm{Im}\,\Gamma \not \subseteq \mathrm{span}\,F$ for any
nontrivial face $F$ of $K$.

It now follows from Lemma~\ref{lemDSR} that $\Gamma V(y)$ is strictly
$K$-quasipositive for each $y \in \mathrm{int}(\mathbb{R}^n_{\geq
0})$. \qquad \qed
\end{pot2}

\section{Applications to chemical reaction networks}

\subsection{Assumptions on reaction kinetics: Condition K}
\label{seckinetics}

We consider chemical reaction networks of the form (\ref{maineqn})
defined on $\mathbb{R}^{m}_{\geq 0}$. At various points above we have
mentioned weak assumptions on the kinetics of reactions, and here we
spell these out. Let $x_i$ refer to the concentration of the $i$th
chemical. Each reaction has a left-hand side and a right-hand side:
let $\mathcal{I}_{j, l}$ be the indices of chemicals occurring on
the left of reaction $j$ and $\mathcal{I}_{j, r}$ be the indices of
the chemicals occurring on the right of reaction $j$. Define
$\Gamma_{ij}$ to be the difference between the number of molecules of
the $i$th chemical occurring on the right of reaction $j$ and on the
left of reaction $j$. Let $v_j$ represent the rate at which reaction
$j$ proceeds to the right. With these definitions we indeed get the
dynamical system (\ref{maineqn}) for the evolution of the
concentrations. We assume that reactions are either reversible or
irreversible, and for an irreversible reaction we make the arbitrary
choice that the substrates occur on the left and products on the
right so the reaction always proceeds with nonnegative velocity to
the right.

We assume:
\begin{enumerate}
\item[(K1)]
$\Gamma_{ij}(\partial v_j/\partial x_i) \leq 0$, and if $\Gamma_{ij} = 0$ then $\partial v_j/\partial x_i = 0$. In other words $V(x) \equiv \left[\frac{\partial v_i(x)}{\partial x_j}\right] \in \mathcal{Q}_0(-\Gamma^\mathrm{T})$. This condition has been discussed before in \cite{banajiSIAM,leenheer} for example, and is satisfied by all reasonable kinetics provided no chemical occurs on both sides of
a reaction.

\item[(K2)]
If reaction $j$ is irreversible then:
\begin{enumerate}
\item[(i)]
$v_j \geq 0$ with $v_j = 0$ if and only if $x_i = 0$ for some
$i \in \mathcal{I}_{j, l}$.
\item[(ii)]
If $x_i > 0$ for all $i \in \mathcal{I}_{j,l}$, then
$\partial v_j/\partial x_i > 0$ for each $i \in
\mathcal{I}_{j,l}$.
\end{enumerate}
\item[(K3)]
If reaction $j$ is reversible then:
\begin{enumerate}
\item[(i)]
If $x_{i} = 0$ for some $i \in \mathcal{I}_{j,l}$ (resp. for
some $i \in \mathcal{I}_{j,r}$) then $v_j \leq 0$ (resp. $v_j
\geq 0$).
\item[(ii)]
If $x_{i} = 0$ for some $i \in \mathcal{I}_{j,l}$ (resp. for some $i \in \mathcal{I}_{j,r}$), then $v_j < 0$ (resp. $v_j > 0$) if and only if $x_{i'} > 0$ for each $i' \in \mathcal{I}_{j,r}$ (resp. for each $i' \in \mathcal{I}_{j,l}$).
\item[(iii)]
If $x_i > 0$ for all $i \in \mathcal{I}_{j,l}$ (resp. for all
$i \in \mathcal{I}_{j,r}$) then, for each $i \in
\mathcal{I}_{j,l}$ (resp. $i \in \mathcal{I}_{j,r}$),
$\partial v_j(x)/\partial x_i
> 0$ (resp. $\partial v_j(x)/\partial x_i < 0$).
\end{enumerate}
\end{enumerate}

We will refer to these very mild and natural assumptions collectively
as {\bf Condition K}.  The following simple result will be useful in
the sequel:
\begin{lem}
\label{lemma:invariance}
Let a system of the form (\ref{maineqn}) satisfy Condition K.  Then
for any $x \in \mathbb{R}^{m}_{\geq 0}$, any $j$ and any $i$ such
that $x_i = 0$ there holds $\Gamma_{ij} v_j(x) \ge 0$. Consequently,
for such a system $\mathbb{R}^{m}_{\geq 0}$ is positively invariant.
\end{lem}
\begin{pf}
Suppose the contrary: then there exists a reaction $j$, some $x \in
\mathbb{R}^{m}_{\geq 0}$ and a chemical $i$ such that $x_i = 0$ and
$\Gamma_{ij} v_j(x) < 0$.  First of all, if
chemical $i$ does not participate in reaction $j$ then $\Gamma_{ij}
= 0$, so this cannot be the case.  If $i$ occurs on the left side of
$j$ and $j$ is irreversible then, by (K2)(i), $v_j(x) = 0$,
consequently $\Gamma_{ij} v_j(x) = 0$.  If $i$ occurs only on the
right of irreversible $j$ then $\Gamma_{ij} > 0$, and, again by
(K2)(i), $\Gamma_{ij} v_j(x) \ge 0$.  Assume now that reaction $j$ is
reversible.  Since, by application of (K3)(i), if chemical $i$
occurs on both sides of reaction $j$ then $x_i=0$ implies $v_j(x) =
0$, this is excluded. Consequently, either (i)
chemical $i$ occurs only on the left of reaction $j$ and $v_j(x) >
0$, or (ii) chemical $i$ occurs only on the right of reaction $j$
and $v_j(x) < 0$. These possibilities are ruled out by assumption (K3)(i). \qquad \qed
\end{pf}

Given a reversible system of reactions satisfying Condition K,
knowledge of the sign pattern of $\Gamma$ allows us to determine the
sign pattern of $V(x)$ at each $x \in \mathbb{R}^{m}_{\geq 0}$; this
in turn allows construction of the digraph at each $x$. (When some
reactions are irreversible, if we have knowledge of how products
affect reaction rates, then again we can determine the zero pattern
of $V(x)$ from $\Gamma$ alone, and hence construct the digraph.) In
particular, if all reactions are reversible then in fact for each $x
\in \mathrm{int}(\mathbb{R}^{m}_{\geq 0})$, $V(x) \in
\mathcal{Q}(-\Gamma^\mathrm{T})$, and $G_{\Gamma, V(x)} = G_{\Gamma,
\Gamma^T}$.

\subsection{Repelling faces of $\mathbb{R}^m_{\geq 0}$}
Let $S \subseteq \{1, \ldots, m\}$ be a subset of $\{1, \ldots, m\}$,
and $S^c = \{1, \ldots, m\}\backslash S$. Define
\[
F_S = \{x \in \mathbb{R}^m\,|\, x_i > 0, i \in S\,\,\,
\mbox{and}\,\,\, x_i = 0, i \in S^c\}.
\]
Following \cite{dubins}, $F_S$ will be referred to as an {\bf
elementary face} of $\mathbb{R}^m_{\geq 0}$. (Elementary faces are
the relative interiors of the closed faces of $\mathbb{R}^m_{\geq
0}$.) If $S$ is a proper, nonempty subset of $\{1, \ldots, m\}$, then
$F_S$ will be termed a {\bf nontrivial} elementary face of
$\mathbb{R}^m_{\geq 0}$. Given a vector field on $\mathbb{R}^m_{\geq
0}$ defining a local semiflow on $\mathbb{R}^m_{\geq 0}$, $F_S$ is
{\bf repelling} if, at each $x \in F_S$ there exists $i \in S^c$ such
that $\dot x_i > 0$. The following is almost obvious:

\begin{lem}
\label{repellingbdy}
Consider a vector field on $\mathbb{R}^m_{\geq 0}$ which defines a
local semiflow $\phi$ on $\mathbb{R}^m_{\geq 0}$, and assume that
some elementary face $F \subseteq \partial \mathbb{R}^m_{\geq 0}$ is
repelling. Then $F$ contains no $\omega$-limit points of $\phi$.
\end{lem}
\begin{pf}
By the definition of a repelling face, each point on $F$ is a start
point of $\phi$: given $x \in F$ since there exists $i$ such that
$x_i = 0$ and $\dot x_i > 0$, there can clearly be no negative
continuation of the trajectory through $x$. But since
$\mathbb{R}^m_{\geq 0}$ is locally compact, no $\omega$-limit set of
$\phi$ can contain a start point of $\phi$ (Theorem 5.6 in
\cite{bhatiahajek}) and so $F$ cannot intersect any $\omega$-limit
set of $\phi$. \qquad \qed
\end{pf}

{\bf Remark.} The following implication of Condition K will be needed
below. Given an elementary face $F_S$, some $y, y' \in F_S$, some $i
\in S^c$, and any $j$ such that $\Gamma_{ij} \not = 0$, then $v_j(y)
> 0$ (resp. $v_j(y) < 0$, resp. $v_j(y) = 0$) implies $v_j(y') > 0$
(resp. $v_j(y') < 0$, resp. $v_j(y') = 0$).

\begin{lem}
\label{tangentrepel}
Consider a system of the form (\ref{maineqn}) defined on
$\mathbb{R}^{m}_{\geq 0}$ and satisfying Condition K. If for some
elementary face $F_S \subseteq \mathbb{R}^m_{\geq 0}$ there exists $x
\in F_S$ where $f(x)$ is not tangent to $F_S$ (i.e., $f(x) \not \in
\mathrm{span}\,F_S$) then $F_S$ is repelling.
\end{lem}
\begin{pf}
If $f(x)$ is not tangent to $F_S$, then there must exist $i \in S^c$
with $\dot x_i \ne 0$.  As $\dot x_i = \sum_j
\Gamma_{ij}v_j(x)$ and each of the summands is, by
Lemma~\ref{lemma:invariance}, nonnegative, it follows that $\dot{x}_i
> 0$.  As a consequence of the remark above, for any $y \in F_s$,
$\dot y_i = \sum_j \Gamma_{ij}v_j(y) > 0$. Thus $F_S$ is repelling.
\qquad \qed
\end{pf}

%Explanation of why a repelling boundary point $z$ is a start point. Let $\dot x = F(x)$. Suppose there exists $y$ and $r > 0$ such that $\phi_r(y) = z$. Then by continuity $\phi_{t_n}(y) \to z$ when $t_n \to r$. By the definition of the vector field
%\[
%\lim_{t_n \to r}\frac{z - \phi_{t_n}(y)}{|r - t_n|} = F(z)
%\]
%and when $x_i = 0$:
%\[
%\lim_{t_n \to r}\frac{\phi_{t_n}(y_i)}{|r - t_n|} = -F_i(z).
%\]
%But clearly if $\phi_{t_n}(y_i) \geq 0$ then the lhs is nonnegative, so the rhs cannot be negative.

\subsection{A family of globally convergent CRNs}

We consider the following class of CRNs with $n$ reactions on $n+1$
chemicals ($n \geq 2$). We refer to the $n$th member of this class
as $\mathcal{R}^{(n)}$ and its stoichiometric matrix as
$\Gamma^{(n)}$. The CRNs are
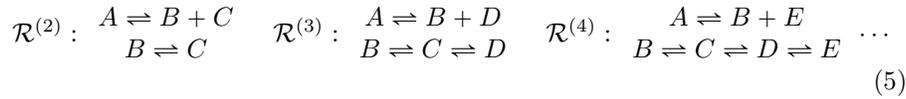
\begin{equation}
\label{eqseq}
\mathcal{R}^{(2)}: \begin{array}{c}A \rightleftharpoons B + C\\ B
\rightleftharpoons C
\end{array}
\quad \mathcal{R}^{(3)}: \begin{array}{c}A \rightleftharpoons B + D\\
B \rightleftharpoons C \rightleftharpoons D
\end{array}
\quad \mathcal{R}^{(4)}: \begin{array}{c}A \rightleftharpoons B + E\\
B \rightleftharpoons C \rightleftharpoons D \rightleftharpoons E
\end{array} \cdots
\end{equation}

Fixing $n$, and choosing a convenient ordering for the reactions the
stoichiometric matrix $\Gamma = \Gamma^{(n)}$ is defined by the
rules:
\begin{enumerate}
\item Column $1$
\[
j=1\,\,:\,\,\left\{\begin{array}{ll}\Gamma_{ij} = -1 & i
=n\\\Gamma_{ij} = 1 & i = n+1\\\Gamma_{ij} = 0 &
\mbox{otherwise}\end{array}\right.
\]

\item Columns $2,\ldots,n-1$
\[
j=2, \ldots,n-1\,\,:\,\,\left\{\begin{array}{ll}\Gamma_{ij} = 1 &
i =n-j+1\\\Gamma_{ij} = -1 & i = n-j+2\\\Gamma_{ij} = 0 &
\mbox{otherwise}\end{array}\right.
\]
\item Column $n$
\[
j=n:\,\,\left\{\begin{array}{ll}\Gamma_{ij} = 1 & i
=1\\\Gamma_{ij} = -1 & i = 2\\\Gamma_{ij} = -1 & i =
n+1\\\Gamma_{ij} = 0 & \mbox{otherwise}\end{array}\right.
\]
\end{enumerate}
$\Gamma^{(n)}$ for $n=2,3,4$ are
\[
\Gamma^{(2)} =
\left(\begin{array}{rr}0&1\\-1&-1\\1&-1\end{array}\right), \quad
\Gamma^{(3)} =
\left(\begin{array}{rrr}0&0&1\\0&1&-1\\-1&-1&0\\
1&0&-1\end{array}\right),
\]
\[
\Gamma^{(4)} =
\left(\begin{array}{rrrr}0&0&0&1\\0&0&1&-1\\0&1&-1&0\\-1&-1&0&0
\\
1&0&0&-1\end{array}\right).
\]

\begin{figure}[h]
\begin{minipage}{\textwidth}
\begin{minipage}{0.3\textwidth}
\begin{center}
\begin{tikzpicture}[domain=0:4,scale=0.4]
\draw[very thin,color=black!0] (-3.5,-3) grid (2.5,3);
\fill (-1.5,1.5) circle (5pt);
\node at  (1.5,1.5) {$B$};

\node at (-3.3,2.4) {$A$}; \draw[<->, line width=0.04cm] (-2.8,2.25)
-- (-1.9, 1.75); %%
\fill (1.5,-1.5) circle (5pt); \node at (-1.5,-1.5) {$C$};

\path (135-15: 2.12cm) coordinate (A1end); \draw[<->, line
width=0.04cm] (A1end)  arc (135-15:45+15:2.12cm);

\path (225-15: 2.12cm) coordinate (A2end); \draw[<->, line
width=0.04cm] (A2end)  arc (225-15:135+15:2.12cm);

\path (315-15: 2.12cm) coordinate (A3end); \draw[<->, line
width=0.04cm] (A3end)  arc (315-15:225+15:2.12cm);

\path (45-15: 2.12cm) coordinate (A4end); \draw[<->, line
width=0.04cm] (A4end)  arc (45-15:-45+15:2.12cm);

\node at (0,0) {$n=2$};
\end{tikzpicture}
\end{center}
\end{minipage}
\begin{minipage}{0.3\textwidth}
\begin{tikzpicture}[domain=0:12,scale=0.48]

%\draw[very thin,color=black!40] (-4,-2) grid (4,2);

\node at (0,0) {$n = 3$};

\path (90: 2cm) coordinate (A2); \path (30: 2cm) coordinate (R2);
\path (-30: 2cm) coordinate (A3); \path (270: 2cm) coordinate (R3);
\path (210: 2cm) coordinate (A1); \path (150: 2cm) coordinate (R1);

\path (-3.8, 1.5) coordinate (B1);

\node at (A1) {$D$}; \node at (A2) {$B$}; \node at (A3) {$C$};

\fill (R1) circle (4pt); \fill (R2) circle (4pt); \fill (R3) circle
(4pt);

\node at (B1) {$A$};

\path (90-15: 2cm) coordinate (A2end); \path (-30-15: 2cm) coordinate
(A3end); \path (210-15: 2cm) coordinate (A1end);

\path (90+15: 2cm) coordinate (A2end1); \path (-30+15: 2cm)
coordinate (A3end1); \path (210+15: 2cm) coordinate (A1end1);

\draw[<->, line width=0.04cm] (A2end)  arc (90-15:90-50:2cm);
\draw[<->, line width=0.04cm] (A3end)  arc (-30-15:-30-50:2cm);
\draw[<->, line width=0.04cm] (A1end)  arc (210-15:210-50:2cm);

\draw[<->, line width=0.04cm] (A2end1)  arc (90+15:90+50:2cm);
\draw[<->, line width=0.04cm] (A3end1)  arc (-30+15:-30+50:2cm);
\draw[<->, line width=0.04cm] (A1end1)  arc (210+15:210+50:2cm);

\draw[<->, line width=0.04cm] (-2.2, 1.1) -- (-3.2, 1.3);

\end{tikzpicture}
\end{minipage}
\begin{minipage}{0.3\textwidth}

\begin{tikzpicture}[domain=-6:6,scale=0.22]

\path (0,0) coordinate (origin);

\node at (0,0) {$n=4$};

\path (45: 6cm) coordinate (R2); \path (-45: 6cm) coordinate (R3);
\path (135: 6cm) coordinate (R1); \path (225: 6cm) coordinate (R4);

\path (90: 6cm) coordinate (A2);

\path (0: 6cm) coordinate (A3); \path (270: 6cm) coordinate (A4);
\path (180: 6cm) coordinate (A1);

\node at (A1) {$E$}; \node at (A2) {$B$}; \node at (A3) {$C$}; \node
at (A4) {$D$};

\fill (R1) circle (8pt); \fill (R2) circle (8pt); \fill (R3) circle
(8pt); \fill (R4) circle (8pt);

\path (90-12: 6cm) coordinate (A2end); \path (0-12: 6cm) coordinate
(A3end); \path (270-12: 6cm) coordinate (A4end); \path (180-12: 6cm)
coordinate (A1end);

\path (90+12: 6cm) coordinate (A2end1); \path (0+12: 6cm) coordinate
(A3end1); \path (270+12: 6cm) coordinate (A4end1); \path (180+12:
6cm) coordinate (A1end1);

\draw[<->, line width=0.04cm] (A2end) arc (90-12:90-35:6cm);
\draw[<->, line width=0.04cm] (A3end) arc (-12:-35:6cm); \draw[<->,
line width=0.04cm] (A4end) arc (270-12:270-35:6cm); \draw[<->, line
width=0.04cm] (A1end) arc (180-12:180-35:6cm);

\draw[<->, line width=0.04cm] (A2end1) arc (90+12:90+35:6cm);
\draw[<->, line width=0.04cm] (A3end1) arc (12:35:6cm); \draw[<->,
line width=0.04cm] (A4end1) arc (270+12:270+35:6cm); \draw[<->, line
width=0.04cm] (A1end1) arc (180+12:180+35:6cm);

\path (-8,5) coordinate (B1); \path (8,5) coordinate (B2); \path
(8,-5) coordinate (B3);

\node at (B1) {$A$};

\draw[<->, line width=0.04cm] (-5.1, 4.5) -- (-6.9,4.8);

\end{tikzpicture}
\end{minipage}
\end{minipage}

\caption{\label{SRsequence} The digraphs associated with reaction
networks in (\ref{eqseq}) for $n = 2,3,4$ under the assumption that
$\left[\frac{\partial v_i(x)}{\partial x_j}\right] \in \mathcal{Q}(-\Gamma^T)$. Vertices associated with rows of
$\Gamma$ (i.e., with particular chemicals) are labelled, while those
associated with columns (i.e., with reactions) are filled circles.
Antiparallel edge-pairs are shown as bidirected edges. As the only
matter of interest is connectivity of the graphs, edge signs and
labels are omitted.}
\end{figure}
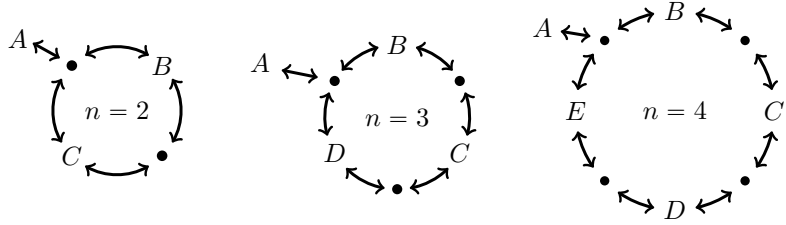

\begin{thm}
\label{appthm1}
Assuming that all reactions are reversible, and the kinetics fulfils
Condition K, the CRN $\mathcal{R}^{(k)}$ ($k \geq 2$) has the
following behaviour: each stoichiometry class contains a unique
equilibrium to which all trajectories are attracted.
\end{thm}
\begin{pf}
The result follows from Theorem~\ref{applicationthm} if all the
conditions in Theorem~\ref{mainthm0} can be shown to be met, and
additionally, $\partial\,\mathbb{R}^{n+1}_{\geq 0}\backslash\{0\}$
contains no limit points of the local semiflow. Note first that for
each $n$, $\Gamma^{(n)}$ is an $(n+1) \times n$ matrix and it is
straightforward that $\mathrm{ker}\,\Gamma^{(n)}$ is trivial. We fix
$n \geq 2$, and for notational simplicity omit the superscript $(n)$
from all objects. Defining $B$ as before, let $\Lambda$ be the $(n+1)
\times 2^{n}$ matrix:
\[
\Lambda = e_{n}\mathbf{1}^T + \Gamma B.
\]
For example, for $n=3$,
\begin{eqnarray*}
\Lambda& = & \left(\begin{array}{c}
0\\0\\1\\0\end{array}\right)
\begin{array}{c}(1,1,1,1,1,1,1,1)\\{}\\{}\\{}\end{array}
+ \\
&& \hspace{3cm}\left(\begin{array}{rrr}
0&0&1\\0&1&-1\\-1&-1&0\\1&0&-1\end{array}\right)
\left(\begin{array}{cccccccc}
0&1&0&1&0&1&0&1\\
0&0&1&1&0&0&1&1\\
0&0&0&0&1&1&1&1\end{array}\right) \\
& = & \left(\begin{array}{rrrrrrrr}
\,\,0&\,\,0&\,\,0&0&1&1&1&1\\0&0&1&1&-1&-1&0&0 \\
1&0&0&-1&1&0&0&-1 \\
0&1&0&1&-1&0&-1&0\end{array}\right).
\end{eqnarray*}

\begin{enumerate}
\item
Note that $\{\Lambda_i\}$ are of the form
\[
e_n + \sum_{i \in S} \Gamma_i
\]
where $S$ can be any subset of $\{1, \ldots, n\}$. We have:
\[
\begin{array}{l}
e_n = \Lambda_1, \\
e_{n+1} = e_n + \Gamma_1 = \Lambda_2, \\
e_{n-k} = e_n + \sum_{i=2}^{k+1}\Gamma_i = \Lambda_{2^k + 1}, \,\,k = 1, \ldots, n-2 \\
e_1 = e_n+ e_{n+1} + \sum_{i=2}^n\Gamma_i = \Lambda_2 +
\Lambda_{2^n-1}.
\end{array}
\]
Each of these facts can quickly be proved by induction. Consequently, there exists a nonnegative matrix $P$ such that $\Lambda P = I$.  For~example, for $n = 3$,
\begin{equation*}
P =
\begin{pmatrix}
0 & 0 & 1 & 0 \\
1 & 0 & 0 & 1 \\
0 & 1 & 0 & 0 \\
0 & 0 & 0 & 0 \\
0 & 0 & 0 & 0 \\
0 & 0 & 0 & 0 \\
1 & 0 & 0 & 0 \\
0 & 0 & 0 & 0
\end{pmatrix}.
\end{equation*}
Thus Condition 1(a) of Theorem~\ref{mainthm0} is met.
\item
It is easy to confirm that the $k$th entry of any sum of columns
of $\Gamma$ has value $-1, 0$ or $1$ for $k \not = n$ and value
$-2, -1$ or $0$ for $k = n$. Consequently, $\Lambda$ is a
$(-1,0,1)$-matrix. Thus Condition 1(b) of Theorem~\ref{mainthm0}
is met.
\item
An immediate consequence of Condition K is that for each $x \in
\mathbb{R}^{n+1}_{\geq 0}$, $V(x) \equiv \left[\partial v_i(x)/
\partial x_j \right] \in \mathcal{Q}_0(-\Gamma^\mathrm{T})$. Thus
Condition 2 of Theorem~\ref{mainthm0} is met.
\item
Using (K3)(iii), the digraphs $G_{\Gamma, V(x)}$ at any
point $x \in \mathrm{int}(\mathbb{R}^{n+1}_{\geq 0})$ can be
inferred. These are illustrated for $n = 2, 3, 4$ in
Figure~\ref{SRsequence}. Clearly they are strongly connected,
each consisting of a single double-cycle with one offshoot. Thus
Condition 3 of Theorem~\ref{mainthm0} is met.
\item
Each nontrivial elementary face of $\mathbb{R}^{n+1}_{\geq 0}$ is
repelling. If this is not so for some nontrivial elementary face
$F_S \subseteq \mathbb{R}^{n+1}_{\geq 0}$, then by
Lemma~\ref{tangentrepel}, $\Gamma v(x)$ is tangent to $F_S$ at
each $x \in F_S$. This in turn implies for each $x \in F_S$ and
each $i \in S^c$, that $\dot x_i = 0$.  Since, by
Lemma~\ref{lemma:invariance}, $i \in S^c$ implies $\Gamma_{ij}
v_j(x) \ge 0$ for each $j$, $\dot x_i = 0$ now implies that
$\Gamma_{ij} v_j(x) = 0$ for each $j$. Thus whenever $\Gamma_{ij}
\not = 0$ then $v_j(x)=0$.  Consider some $j$ such that
$\Gamma_{ij} \not = 0$. Note that no chemical occurs on more
than one side of this reaction system and that all reactions are
reversible. Suppose that chemical $i$ appears only on the right
of reaction $j$ so that $\Gamma_{ij} > 0$.  $v_j(x) = 0$ then
implies, via (K3)(ii), that there exists $i' \in \mathcal{I}_{j,
l}\cap S^c$, and so $\Gamma_{i'j} < 0$. Suppose, on the other
hand that chemical $i$ appears only on the left of reaction $j$
so that $\Gamma_{ij} < 0$. In this case, via (K3)(ii), there
exists $i' \in \mathcal{I}_{j, r}\cap S^c$ and so $\Gamma_{i'j} >
0$. Thus, for a system of reversible reactions such as this, with
no chemical appearing on more than one side of any reaction, the
condition that $\Gamma v(x)$ is tangent to $F_S$ at each $x \in
F_S$ has a simple interpretation in terms of $\Gamma$: each
column of the submatrix $\Gamma^{S^c}$ of $\Gamma$ with only rows
indexed by $S^c$ must contain either both a positive and a
negative entry, or only zeros. That this does not occur in
$\Gamma^{(n)}$ for any nonempty $S \subsetneq \{1, \ldots, n+1\}$
can be checked. If $S$ includes $1$, it must include $2$ or
$n+1$; if $S$ includes $k$ where $2\leq k\leq n$, then it must
include $k+1$ and $k-1$; if $S$ includes $n+1$, then it must
include $1$; thus if $S$ is nonempty then $S = \{1, \ldots,
n+1\}$.  In the terminology of \cite{angelipetrinet}, the systems
have no ``siphons''. In physical terms, setting any proper subset
of reactant concentrations to zero, there always remains some
active reaction producing one of the absent reactants with
positive rate. Consequently, by Lemma~\ref{repellingbdy}, no
nontrivial elementary face of $\mathbb{R}^{n+1}_{\geq 0}$
contains any $\omega$-limit points of the associated semiflow
$\phi$. As $\partial \mathbb{R}^{n+1}_{\geq 0}\backslash\{0\}$ is
precisely the union of these nontrivial elementary faces, it
contains no $\omega$-limit points of the semiflow.
\end{enumerate}
\qquad \qed
\end{pf}

\subsection{Another example}

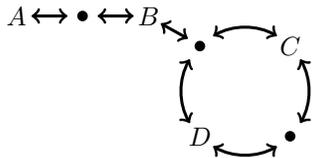
\begin{figure}[h]
\begin{center}
\begin{tikzpicture}[domain=0:4,scale=0.4]
\draw[very thin,color=black!0] (-3.5,-3) grid (2.5,3);
\fill (-1.5,1.5) circle (5pt);
\node at (1.5,1.5) {$C$};

\node at (-3.2,2.5) {$B$}; \draw[<->, line width=0.04cm] (-2.8,2.25)
-- (-1.9, 1.75); %%
\fill (1.5,-1.5) circle (5pt); \node at (-1.5,-1.5) {$D$};

\path (135-15: 2.12cm) coordinate (A1end); \draw[<->, line
width=0.04cm] (A1end)  arc (135-15:45+15:2.12cm);

\path (225-15: 2.12cm) coordinate (A2end); \draw[<->, line
width=0.04cm] (A2end)  arc (225-15:135+15:2.12cm);

\path (315-15: 2.12cm) coordinate (A3end); \draw[<->, line
width=0.04cm] (A3end)  arc (315-15:225+15:2.12cm);

\path (45-15: 2.12cm) coordinate (A4end); \draw[<->, line
width=0.04cm] (A4end)  arc (45-15:-45+15:2.12cm);

\fill (-5.4,2.5) circle (5pt); \node at (-7.6,2.5) {$A$}; \draw[<->,
line width=0.04cm] (-3.7,2.5) -- (-4.9, 2.5); \draw[<->, line
width=0.04cm] (-5.9,2.5) -- (-7.1,2.5);

\end{tikzpicture}
\end{center}
\caption{\label{SRlone} The digraph $G_{\Gamma, V}$ for the reaction
system in (\ref{eqtail}). Conventions are as in
Figure~\ref{SRsequence}.}
\end{figure}

The family of CRNs presented in the previous section are not the only
CRNs preserving cubic orders. As another example of a CRN which preserves a cubic order, consider
the following:
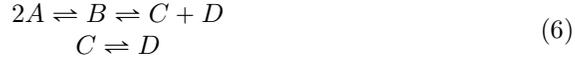
\begin{equation}
\label{eqtail}
\begin{array}{c}
2A \rightleftharpoons B \rightleftharpoons C + D \\
C \rightleftharpoons D\end{array}
\end{equation}
giving
\[
\Gamma = \left(\begin{array}{rrr}
0&0&2 \\
0&1&-1 \\
-1&-1&0 \\
1&-1&0
\end{array}\right).
\]
$\mathrm{ker}\,\Gamma = \{0\}$ and the positive vector $(1,2,1,1)^T$
lies in $\mathrm{ker}\,\Gamma^\mathrm{T}$. Defining $B$ in the usual
way,
\[
\Lambda \equiv e_3\mathbf{1}^\mathrm{T} + \Gamma B =
\left(\begin{array}{rrrrrrrr}
0&0&0&0&2&2&2&2 \\
0&0&1&1&-1&-1&0&0 \\
1&0&0&-1&1&0&0&-1 \\
0&1&-1&0&0&1&-1&0\end{array}\right).
\]
Left-multiplying $\Lambda$ by the diagonal matrix
$\mathrm{diag}\{1/2, 1, 1, 1\}$ gives a $(-1, 0, 1)$-matrix. The
matrix
\[
P =
\left(\begin{array}{cccc} 1/2&0&1&0 \\
0&1&0&1 \\
0&1&0&0 \\
0&0&0&0 \\
0&0&0&0 \\
0&0&0&0 \\
0&0&0&0 \\
1/2&0&0&0
\end{array}\right)
\]
can be confirmed to be a nonnegative right-inverse for $\Lambda$. Consequently, assumption 1 in Theorem~\ref{mainthm0} is satisfied. Assuming that the kinetics fulfils (K3)(iii), the digraph in $\mathrm{int}(\mathbb{R}^4_{\geq 0})$ is shown in Figure~\ref{SRlone} and is clearly strongly connected.

Thus all the conditions of Theorem~\ref{mainthm0} are fulfilled.
Further, it can easily be checked, using the approach in part 5 of
the proof of Theorem~\ref{appthm1}, that each nontrivial elementary
face of $\mathbb{R}^{n+1}_{\geq 0}$ is repelling and so by
Lemma~\ref{repellingbdy}, $\partial \mathbb{R}^{n+1}_{\geq
0}\backslash\{0\}$ contains no $\omega$-limit points of the
associated semiflow $\phi$. Thus by Theorem~\ref{applicationthm} each
stoichiometry class contains a unique equilibrium to which all
trajectories are attracted. Clearly with the exception of the
equilibrium at $0$, these globally attracting equilibria lie in
$\mathrm{int}(\mathbb{R}^{n+1}_{\geq 0})$.

We remark that in the above examples, the assumption that all
reactions are reversible can be relaxed provided that this is done in
such a way that the digraph in $\mathrm{int}(\mathbb{R}^{n+1}_{\geq
0})$ remains strongly connected, and such that nontrivial elementary 
faces of $\mathbb{R}^{n+1}_{\geq 0}$ remain repelling.

\section{Concluding remarks}

Systems of differential equations arising in chemistry and displaying
strong convergence properties have been presented. The examples are
only intended to illustrate the applicability of the main theorems
(Theorems~\ref{mainthm0}~and~\ref{applicationthm}): characterising
fully the set of CRNs for which global convergence can be proved via
these theorems remains a task for the future.

In addition, a number of variants on the same theme remain to be
explored. While this work has focussed on systems preserving
so-called cubic orders, these are by no means the only partial orders
which may be preserved by a system of chemical reactions as
demonstrated in \cite{banajidynsys}. Closely related results can be
derived where the geometrical construction in this paper is
restricted to linear invariant subspaces: this approach throws up
further examples of chemical reaction networks with strong global
convergence properties. These themes will be explored in forthcoming
papers.

\section*{Acknowledgements}

The work of M. Banaji was supported by EPSRC grant EP/J008826/1 ``Stability and order preservation in chemical reaction networks''. The work of J. Mierczy\'nski was supported by research grant MENII N N201 394537 (2009--2012), Poland.

%\section*{References}

\bibliographystyle{unsrt}

\begin{thebibliography}{10}

\bibitem{hornjackson}
F.~Horn and R.~Jackson.
\newblock General mass action kinetics.
\newblock {\em Arch. Ration. Mech. Anal.}, 47(2):81--116, 1972.

\bibitem{feinberg}
M.~Feinberg.
\newblock Chemical reaction network structure and the stability of complex
  isothermal reactors - {I}. {T}he deficiency zero and deficiency one theorems.
\newblock {\em Chem. Eng. Sci.}, 42(10):2229--2268, 1987.

\bibitem{halsmith}
H.~Smith.
\newblock {\em Monotone Dynamical Systems: An introduction to the theory of
  competitive and cooperative systems}.
\newblock American Mathematical Society, 1995.

\bibitem{hirschsmith}
M.W. Hirsch and H.~Smith.
\newblock {\em Handbook of Differential Equations: Ordinary Differential
  Equations, {Vol II}}, chapter Monotone Dynamical Systems, pages 239--357.
\newblock Elsevier B. V., Amsterdam, 2005.

\bibitem{mierczynski}
J.~Mierczy\'nski.
\newblock Strictly cooperative systems with a first integral.
\newblock {\em SIAM J. Math. Anal.}, 18(3):642--646, 1987.

\bibitem{banajiangeli}
M.~Banaji and D.~Angeli.
\newblock Convergence in strongly monotone systems with an increasing first
  integral.
\newblock {\em SIAM J. Math. Anal.}, 42(1):334--353, 2010.

\bibitem{banajiangelierratum}
M.~Banaji and D.~Angeli.
\newblock Addendum to ``{C}onvergence in strongly monotone systems with an
  increasing first integral''.
\newblock {\em SIAM J. Math. Anal.}, 44(1):536--537, 2012.

\bibitem{berman}
A.~Berman and R.~Plemmons.
\newblock {\em Nonnegative matrices in the mathematical sciences}.
\newblock Academic Press, New York, 1979.

\bibitem{nikaido}
H.~Nikaido.
\newblock {\em Convex structures and economic theory}.
\newblock Academic Press, 1968.

\bibitem{websterconvexity}
R.~Webster.
\newblock {\em Convexity}.
\newblock Oxford University Press, 1994.

\bibitem{brualdi}
R.~A. Brualdi and B.~L. Shader.
\newblock {\em Matrices of sign-solvable linear systems}.
\newblock Number 116 in Cambridge tracts in mathematics. Cambridge University
  Press, 1995.

\bibitem{banajicraciun2}
M.~Banaji and G.~Craciun.
\newblock Graph-theoretic approaches to injectivity and multiple equilibria in
  systems of interacting elements.
\newblock {\em Commun. Math. Sci.}, 7(4):867--900, 2009.

\bibitem{banajiburbanks}
M.~Banaji and A.~Burbanks.
\newblock A graph-theoretic condition for irreducibility of a set of cone
  preserving matrices.
\newblock {\em preprint available at {\tt http://arxiv.org/abs/1112.1653}}.

\bibitem{banajiSIAM}
M.~Banaji, P.~Donnell, and S.~Baigent.
\newblock {$P$} matrix properties, injectivity and stability in chemical
  reaction systems.
\newblock {\em SIAM J. Appl. Math.}, 67(6):1523--1547, 2007.

\bibitem{leenheer}
P.~{De Leenheer}, D.~Angeli, and E.D. Sontag.
\newblock Monotone chemical reaction networks.
\newblock {\em J. Math. Chem.}, 41(3):295--314, 2007.

\bibitem{dubins}
L.~E. Dubins.
\newblock On extreme points of convex sets.
\newblock {\em J. Math. Anal. Appl.}, 5:237--244, 1962.

\bibitem{bhatiahajek}
N.~P. Bhatia and O.~Hajek.
\newblock {\em Local semi-dynamical systems}.
\newblock Springer-Verlag, 1969.

\bibitem{angelipetrinet}
D.~Angeli, P.~{De Leenheer}, and E.~D. Sontag.
\newblock A {P}etri net approach to the study of persistence in chemical
  reaction networks.
\newblock {\em Math. Biosci.}, 210:598--618, 2007.

\bibitem{banajidynsys}
M.~Banaji.
\newblock Monotonicity in chemical reaction systems.
\newblock {\em Dyn Syst}, 24(1):1--30, 2009.

\end{thebibliography}

\end{document}